\documentclass[oneside,11pt]{article}
\usepackage{amssymb,amsthm,amsmath,marvosym,amsfonts,fontenc}
\usepackage[usenames,dvipsnames]{color}
\usepackage{enumerate,marginnote}
\usepackage[dvips]{graphicx}

\usepackage{mathtools}
\usepackage{verbatim} 
\usepackage{authblk} 

\usepackage{caption}
\usepackage{fullpage} 
\usepackage{epsfig,psfrag}

\numberwithin{equation}{section}
\numberwithin{figure}{section}
\newtheorem{theorem}{Theorem}[section]
\newtheorem{lem}[theorem]{Lemma}
\newtheorem{cl}[theorem]{Claim}

\newtheorem{prop}[theorem]{Proposition}

\newtheorem{definition}[theorem]{Definition}

\usepackage{authblk}

\begin{document}

\title{A median-type condition for graph tiling}
%

\author{Diana Piguet}
\affil{The Czech Academy of Sciences, Institute of Computer Science}
\author{Maria Saumell\thanks{
		Piguet and Saumell were supported by the Czech Science Foundation, grant number GJ16-07822Y and by institutional support~RVO:67985807. 
		Saumell was also supported by Project LO1506 of the Czech Ministry of Education, Youth and Sports. An extended abstract of this result appears in the  proceedings of EuroComb 2017, \emph{Electronic Notes in Discrete Mathematics
			Volume 61, August 2017, Pages 979-985.} E-mail adresses: \texttt{piguet/saumell@cs.cas.cz}}}
\affil{The Czech Academy of Sciences, Institute of Computer Science\\
	\&
\\Department of Theoretical Computer Science
Faculty of Information Technology
Czech Technical University in Prague}

\date{}

\maketitle

\begin{abstract}
	Koml\'os [Tiling Tur\'an theorems, Combinatorica, 20,{\bf 2} (2000), 203--218] determined the asymptotically
	optimal minimum degree condition for covering a given proportion of vertices of a host
	graph by vertex-disjoint copies of a fixed  graph. We show that the minimum degree condition can be relaxed in the sense that we require only a given fraction of vertices to have the prescribed degree.  
\end{abstract}

\section{Introduction and Results}\label{sec:intro}

A classical question in extremal graph theory is to study which density condition forces a graph to contain a certain fixed subgraph. Let us cite
 Tur\'an theorem~\cite{Turan1941}, which  gives the average degree forcing a clique, or  the Erd\H os--Stone theorem~\cite{Erdos1946}, which essentially determines the average degree condition guaranteeing the containment of a fixed non-bipartite graph~$H$.
This approach can be generalised by seeking the density condition forcing the containment of several copies of~$H$. One might consider edge-disjoint copies of~$H$, a so-called \emph{packing} of~$H$, or  vertex-disjoint copies of~$H$, known as an \emph{$H$-tiling}. The latter concept generalises the one of a matching, where we take~$H$ to be an edge, and this is the direction we follow in this paper.

\begin{definition}[Tiling, perfect tiling, tiling number]\label{def:tiling}
	Let~$H$ be a fixed graph. We say that a graph~$G$ has an \emph{$H$-tilling of size $m$}, if there are $m$ vertex-disjoint copies of~$H$ in~$G$.
	An $H$-tiling is \emph{perfect} if it has size $\frac{|V(G)|}{|V(H)|}$.
	The size of the largest $H$-tiling in~$G$ is called the \emph{$H$-tiling number} (in~$G$) (or just tiling number, when $H$ is clear from the context). 
\end{definition}

The concept is not new and researchers have obtain many results in this direction.
The Hajnal-Szemer\'edi theorem on equitable colouring~\cite{Hajnal1970} (formulated for the complement of the host graph) determines the minimum degree condition guaranteeing a perfect clique tiling. One can deduce from their result the minimum degree condition needed to force a partial clique-tiling of any given size.
For the corresponding question of determining the average degree condition forcing a clique tiling of a given size, the current state of knowledge  is rather poor.
The only two known cases for which the average degree condition has been determined correspond to $K_2$-tilings and $K_3$-tilings, given by Erd\H os and Gallai~\cite{Erdos1959} and  Allen et al.~\cite{ABHP:DensityCorHaj}, respectively.
The case of tiling cliques of higher order is wide open and no conjecture has even been formulated.

Generalising the concept of clique tilings to the one of tilings with an arbitrary graph~$H$, Koml\' os~\cite{Komlos2000} extended the Hajnal-Szemer\'edi theorem.

\begin{theorem}[Koml\'os tiling theorem] \label{thm:Komlos}Let $x\in (0,\frac{1}{|V(H)|})$ and $r\ge 2$.
If~$H$ is an $r$-colourable graph with colour class sizes $\ell_1\ge \ldots\ge \ell_r>0$  and $G$ is an $n$-vertex graph with minimum degree at least $(r-2+x\ell_r)\frac{n}{r-1}$, then $G$ contains an $H$-tiling of size at least~$(x-o(1))n$.
\end{theorem}

As~$H$ may have many $r$-colourings, in order to apply Koml\'os' theorem with the weakest hypothesis, we need to fix an $r$-colouring of~$H$ which minimises the size of its smallest colour class. For the corresponding value of $\ell_r$ and for each $x$, Koml\'os constructs graphs with minimum degree $(r-2+x\ell_r)\frac{n}{r-1}$ and $H$-tiling number~$xn$. This shows that his result is asymptotically optimal. Recently, Hladk\'y, Hu, and Piguet~\cite{HlHuPi:Komlos} proved stability of this result.
Koml\'os tiling theorem was complemented by K\"uhn and Osthus~\cite{Kuehn2009}, who established the optimal minimum degree condition forcing a perfect $H$-tiling up to an additive constant. Considering average degree rather than minimum degree, Grosu and Hladk\'y~\cite{Grosu2012} extended  the Erd\H os--Gallai theorem, by asymptotically establishing the average degree condition forcing the containment of an $H$-tiling, for a fixed bipartite graph~$H$.

A finer approach to extremal problems is to take into account more information encoded in the degree sequence of the host graph. Let us give an example in the area of tree containment. It is trivial to see that a minimum degree of~$k$ ensures a copy of any tree of size~$k$. However,  Loebl, Koml\'os, and S\'os conjecture that only half of the vertices need to have this degree to guarantee the same assertion. Taking this perspective in the area of tilings, Treglown asymptotically determined an optimal degree  sequence condition forcing a perfect $H$-tiling~\cite{MR3471843} and very recently Hyde, Liu and Treglown~\cite{Hyde_Liu_Treglown_degree_sequence} asymptotically determined the optimal degree sequence for a partial tiling.
In this paper,  we inquire
what portion of vertices need to meet the  degree bound in Koml\'os tiling theorem in order to guarantee an $H$-tiling of a given size without any requirement on the other vertices. Our main result is the following.

\begin{theorem}\label{thm:result}
	Let $H$ be an $r$-colourable graph with colour class sizes $\ell_1\ge \ldots\ge \ell_r>0$ and let $x\in (0,\frac{1}{|V(H)|})$. Then for any $\eta>0$ there is an $n_0\in \mathbb N$ such that for any $n\ge n_0$ and for $\delta:= (r-2+x\ell_r)\frac{n}{r-1}$ the following holds.
	
	Any $n$-vertex graph with at least $(1+\eta)(r-2+x|V(H)|)\frac{n}{r-1}$ vertices of degree at least~$(1+\eta)\delta$ contains an $H$-tiling of size at least~$xn$.

\end{theorem}

Theorem~\ref{thm:result} strengthens Koml\'os tiling theorem: the degree bound $\delta$ is the same but we do not require all the vertices of the host graph to meet this bound, but rather only roughly a $f=\frac{r-2+x|V(H)|}{r-1}$ fraction of them. Note that as $x$ ranges from~$0$ to $\frac{1}{|V(H)|}$, $f$ ranges from $\frac{r-2}{r-1}$ to~1. 

\begin{figure}[th] \centering
	\includegraphics[scale=1.00]{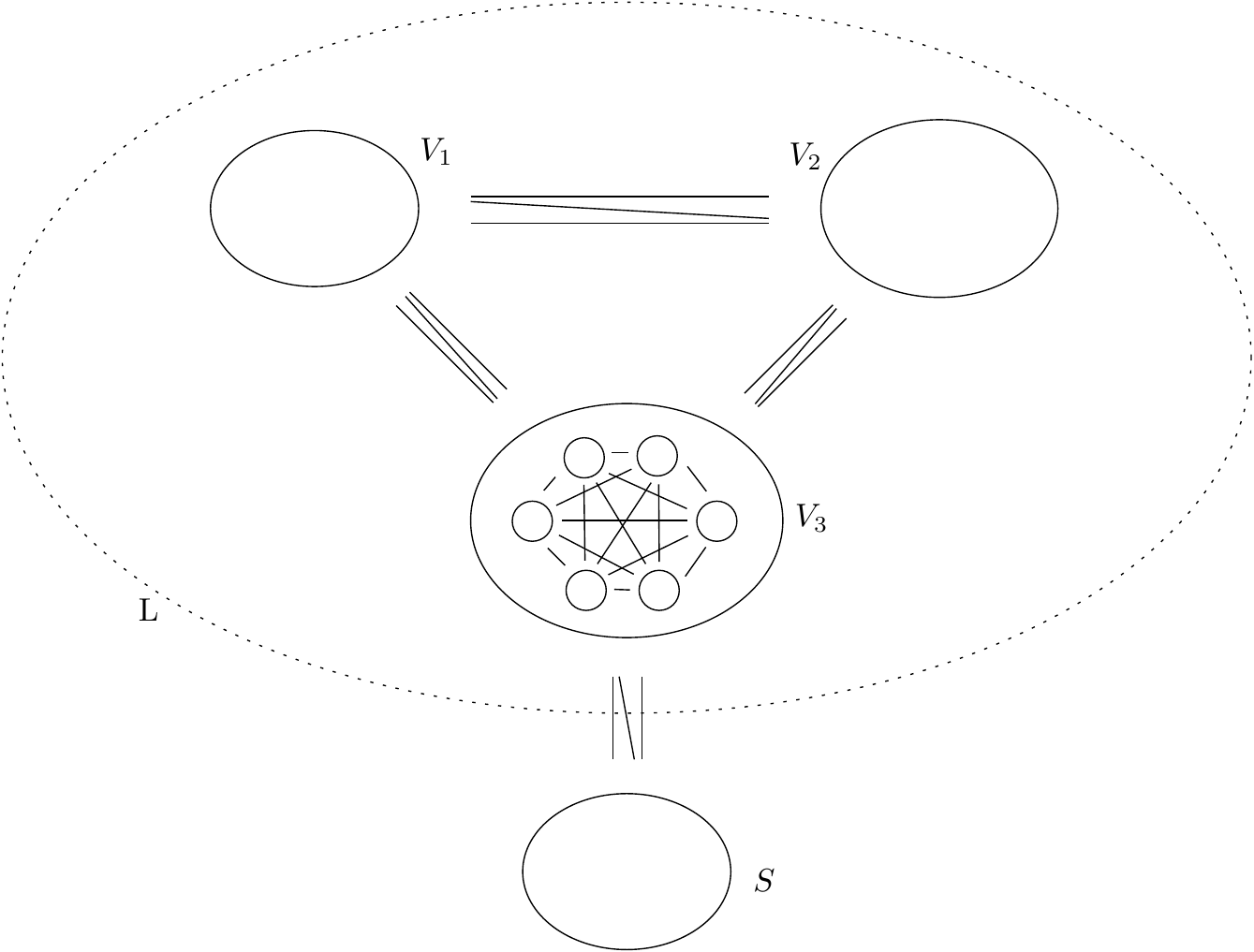}
	\caption{The extremal graph $G(x,|V(H)|,r,\ell_r,n)$. The sizes of the sets are as follows: $|V_1|=x\ell_r n$, $|V_2|=x(|V(H)|-\ell_r)\frac{n}{r-1}$, $|V_3|=(r-2)(1-x\ell_r)\frac{n}{r-1}$, and $|S|=n-|L|=(r-2-x|V(H)|)\frac{n}{r-1}$.}
	\label{fig:extremal}
\end{figure}

When $\ell_r$ is chosen according to the colouring minimising the size of the smallest colour class, Theorem~\ref{thm:result} is asymptotically optimal for each value of $x$. To show this, we construct an $n$-vertex graph $G=G(x,|V(H)|,r,\ell_r,n)$ with  $(r-2+x|V(H)|)\frac{n}{r-1}$ vertices of degree  $\delta:= (r-2+x\ell_r)\frac{n}{r-1}$ that does not contain an $H$-tiling of size more than $xn$. The vertex set of $G$ is partitioned into four sets, $|V_1|=x\ell_r n$, $|V_2|=x(|V(H)|-\ell_r)\frac{n}{r-1}$, $|V_3|=(r-2)(1-x\ell_r)\frac{n}{r-1}$, and $|S|=(1-x|V(H)|)\frac{n}{r-1}$. The sets~$V_1$,~$V_2$ and~$S$ are independent, and the set~$V_3$ induces a balanced complete  $(r-2)$-partite graph. 
There is a complete bipartite graph between~$V_3$ and $V(G)\setminus V_3$ and a complete bipartite graph between~$V_1$ and~$V_2$, and no other edge. The graph $G(x,|V(H)|,r,\ell_r,n)$ is depicted in Figure~\ref{fig:extremal}. Note that the  set $L=V_1\cup V_2\cup V_3$ has size $|L|=(r-2+x|V(H)|)\frac{n}{r-1}$ and its vertices meet our degree assumption. Indeed, vertices in $V_1$ have degree $|L|-|V_1|\ge (r-2+x\ell_r)\frac{n}{r-1}$, as $|V(H)|\ge r\ell_r$, vertices in $V_2$ have degree $|L|-|V_2|=(r-2+x\ell_r)\frac{n}{r-1}$, and vertices~$V_3$ have degree $n-\frac{|V_3|}{r-2}=(r-1-\frac{(r-2)(1-x\ell_r)}{r-2})\frac{n}{r-1}=(r-2+x\ell_r)\frac{n}{r-1}$.
We claim that each copy of~$H$ in~$G$ must have at least~$\ell_r$ vertices in~$V_1$. Suppose it is not the case. We can then recolour~$H$ so that the smallest colour class has less than~$\ell_r$ vertices,  as follows. Give the vertices embedded in~$V_1$ colour~$1$, the ones embedded in~$V_2\cup S$ colour~$2$, and distribute the $r-2$ other colours to vertices embedded in~$V_3$. In this colouring, the smallest colour class has at most as many vertices as colour~$1$ has, which is less than~$\ell_r$. 
Since $|V_1|=x\ell_r n$, we conclude that there cannot be an $H$-tiling of size more than~$xn$.

The proof of Theorem~\ref{thm:result} is postponed to Section~\ref{sec:proof_result} and it relies on the regularity method and on the generalisation of the LP-duality between fractional matching and fractional cover. Such a generalisation of the  LP-duality has already been used in~\cite {Martin2013} for the particular case of \emph{fractional clique-tiling}.  However, we think that this method deserves more attention and is exposed here in the more general form for \emph{fractional $H$-tiling}\footnote{For a precise definition of fractional tiling, see Section~\ref{ssec:duality}.}.
Similar use of LP-duality for graph-tiling was used in~\cite{DoHl:Polytons,HlHuPi:Komlos,HlHuPi:TilingsInGraphons} but in the context of graphons.

\section{Notation and Preliminaries}\label{sec:prelim}
\subsection{The regularity lemma and the embedding lemma}\label{ssec:Regularity}

Given two disjoint subsets $X,Y\subset V(H)$, we denote by $e(X,Y)$ the number of edges between vertices in $X$ and  vertices in $Y$, that is,
\[e(X,Y):=\left|\left\{(x,y)\in E(H),\ x\in X\ \mbox{and}\ y\in Y\right\}\right|\,.\]

\begin{definition}[Regular pairs, bipartite density]
For a given $\varepsilon>0$, a pair $(U,W)$ of disjoint
sets $U,W\subset V(H)$ 
is called an \emph{$\varepsilon$-regular
	pair} if  for every
$U'\subset U$, $W'\subset W$ with $|U'|\ge \varepsilon |U|$, $|W'|\ge
\varepsilon |W|$ we have that $|d(U,W)-d(U',W')|<\varepsilon$, where the \emph{(bipartite) density} $d$ between two disjoint sets $X$ and $Y$ is defined as $d(X,Y):= \frac{e(X,Y)}{|X|\cdot |Y|}$. If the pair $(U,W)$ is not $\varepsilon$-regular,
then we call it \emph{$\varepsilon$-irregular}. 
\end{definition}

\begin{lem}[Slicing lemma, Fact~1.5 in~\cite{Komlos1996a}]\label{lem:subdivide}
	Let $G$ be a graph and $(X,Y)$ be an $\varepsilon$-regular pair of density $d$ in $G$. Let $X'\subseteq X$ and $Y'\subseteq Y$ such that $|X'|\geq \alpha |X|$ and $|Y'|\geq \alpha |Y|$. Then, $(X',Y')$ is an $\varepsilon'$-regular pair of density at least $d-\varepsilon$, where $\varepsilon'=\max(\varepsilon/\alpha,2\varepsilon)$.
\end{lem}

For completeness, let us next state the fundamental regularity lemma.

\begin{lem}[Szemer\'edi's regularity lemma, Theorem~7.4.1 in~\cite{Diestel2016}]\label{lem:RL}
	For all $\varepsilon_R>0$ and $\ell\in\mathbb{N}$ there exist $n_R,M\in\mathbb{N}$
	such that for every $n\ge n_R$ and every $n$-vertex graph~$G$  there exists
	a partition $U_0,U_1,\ldots,U_p$ of $V(G)$, $\ell<p<M$, with the following properties:
	\begin{enumerate}[1)]
		\item For every $i,j\in [p]$ we have $|U_i|=|U_j|:=m$, and  $|U_0|<\varepsilon_R n$.
		\item All but at most
		$\varepsilon_R p^2$ pairs $(U_i,U_j)$, $i,j\in [p]$, $i\neq j$, are $\varepsilon_R$-regular.
	\end{enumerate}
\end{lem} %

The partition $U_0,U_1,\ldots,U_p$ of $V(G)$ allows to define a \emph{cluster graph} $\mathbf G$ (known as well as \emph{reduced graph}) with parameters $\varepsilon_R$, $m$ and $d$ as follows: $V(\mathbf G):=\{U_1,\ldots,U_p\}$ and $U_iU_j\in E(\mathbf G)$ if and only if $(U_i,U_j)$ is an $\varepsilon_R$-regular pair such that $d
(U_i,U_j)\ge d$.

Given a cluster graph $\mathbf G$, let \emph{the blow-up graph} $\mathbf G_s$ denote the graph defined as follows: Every vertex~$v$ of $\mathbf G$ is replaced by a set of $s$ vertices, which we call the \emph{clone set} of~$v$, and every edge by a complete bipartite graph between the corresponding clone-sets.

There is a natural association between the clusters of $G$, the vertices of $\mathbf G$, and the clone-sets of~$\mathbf G_s$: The cluster $U_i$ of $G$ is associated to the vertex $U_i$ of~$\mathbf G$ and, in $\mathbf G_s$, this vertex is replaced by the clone-set of~$U_i$, which we denote by $U_i^{(s)}$

\begin{lem}[Lemma 7.5.2 in~\cite{Diestel2016} (Embedding lemma)]\label{lem:emb}
	 For any $d\in \left( 0,1\right]$ and $\Delta \geq 1$, there exists an $\varepsilon_0> 0$ such that the following holds. Let $\mathbf G$ be the cluster graph of a graph~$G$ with parameters $\varepsilon_R>0$ for the regularity, $m\in \mathbb N$ for the size of its clusters and $d>0$ for the minimal density of the regular pairs. If $\Delta(H)\leq \Delta$,  $\varepsilon_R\leq \varepsilon_0$ and $m\geq 2s/d^\Delta$, then 
	\[H\subseteq \mathbf G_s\ \Rightarrow\ H\subseteq G \,.\]
\end{lem}	

From the proof of Lemma~\ref{lem:emb} given in~\cite{Diestel2016}, one can infer the following slightly stronger claim: If the hypothesis of the lemma are fulfilled and $H\subseteq \mathbf G_s$, then there exists an embedding of~$H$ in~$G$ such that, if a vertex of $H\subseteq \mathbf G_{s}$ is contained in clone-set~$U_j^{(s)}$ of~$\mathbf G_{s}$, then this vertex is contained in the cluster~$U_j$ of~$G$.

\subsection{Fractional hom$_H$-tiling, fractional hom$_H$-cover, and LP-duality}\label{ssec:duality}

From Lemma~\ref{lem:emb} it is clear that a copy of $H$ in a cluster graph~$\mathbf G$ implies a copy of~$H$ in the original graph~$G$; actually using Lemma~\ref{lem:emb} it can be shown that it generates up to $(1-o(1))\frac{|V(G)|}{|V(\mathbf G)|}$ many disjoint copies of~$H$. Going through all (not-necessarily disjoint) copies of~$H$ in~$\mathbf G$, we could use the notion of fractional tiling to find a large~$H$-tiling in~$G$.
A \emph{fractional tiling} of~$H$ in~$G$ is a function giving weights in $[0,1]$ to copies of~$H$ in~$G$ such that for each vertex $v\in V(G)$ the sum of the weights of all copies containing~$v$ is at most~$1$. Hence an $H$-tiling in~$G$ (Definition~\ref{def:tiling}) is an \emph{integral} tiling, i.e., a fractional tiling whose weights are either~$0$ or~$1$. 

But not only copies of~$H$ in~$\mathbf G$ may generate copies of~$H$ in~$G$. For example, if there is a triangle in~$\mathbf G$, we know there are many copies of~$C_5$ in~$G$. Therefore, instead of looking for copies of~$H$ in~$\mathbf G$, we shall seek copies of~$H'$, where~$H'$ is some homomorphic image of~$H$. 
A \emph{homomorphism} $h:H\rightarrow G$ is a mapping from $V(H)$ to $V(G)$ such that 
$uv\in E(H)$ implies $h(u)h(v)\in E(G)$.
Hence, we shall generalise the concept of fractional tiling of~$H$ in~$G$, by enriching it by the homomorphism images of~$H$.

Set 

\[G^H:=\{h\:|\: h:H\rightarrow G\mbox{ is a homomorphism}\}\;.\]

\begin{definition}[Fractional hom$_H$-tiling]\label{def:f H-tiling}
A function $f:\{h(H)\subseteq G\::\: h\in G^H\}\rightarrow [0,1]$ is a fractional hom$_H$-tiling in $G$ of size $\kappa\in \mathbb R^+$ if it satisfies the following two properties:
\begin{enumerate}
\item \label{it:1}For any vertex $v\in V(G)$, we have $\sum f(h(H))\cdot |h^{-1}(v)|\le 1$, where the sum runs over all homomorphisms $h\in G^H$.
\item \label{it:2}We have that $\sum f(h(H))=\kappa$, where the sum runs through all homomorphisms $h\in G^H$.
\end{enumerate}
\end{definition}

So in a fractional hom$_H$-tiling we not only assign weights to isomorphic copies of~$H$ in~$G$, but to homomorphic copies of~$H$, as well.

In order to prove Theorem~\ref{thm:result}, we shall use LP-duality.
So, similarly as above,  we need to generalise the notion of fractional $H$-cover in a graph~$G$ to consider also homomorphic copies of~$H$, which will be the dual notion of fractional hom$_H$-tiling (Definition~\ref{def:f H-cover}). We shall prove in Proposition~\ref{prop:hom-duality} that those two notions are indeed dual.

\begin{definition}[Fractional hom$_H$-cover]\label{def:f H-cover}
A function $c:V(G)\rightarrow [0,1]$ is a  fractional hom$_H$-cover in $G$ of size $\kappa\in \mathbb R^+$ if it satisfies the following two properties:
\begin{enumerate}
\item [(i)] \label{it:i} For any subgraph $H'\subseteq G$ and any homomorphism $h:H\rightarrow H'$, we have $\sum c(v)|h^{-1}(v)|\ge 1$, where the sum runs through all vertices $v\in V(H')$.
\item [(ii)] \label{it:ii} We have $\sum c(v)=\kappa$, where the sum runs through all vertices $v\in V(G)$. 
\end{enumerate}
\end{definition}

\begin{definition}[Fractional tiling/cover number]
The fractional hom$_H$-tiling number of a graph~$G$ is the maximum of the sizes of all its fractional hom$_H$-tilings.

The fractional hom$_H$-cover number of a graph~$G$ is the minimum of the sizes of all its fractional hom$_H$-covers.
\end{definition}

\begin{prop}[LP-duality for hom$_H$]\label{prop:hom-duality}
For any graph~$G$, its fractional hom$_H$-tiling number equals its fractional hom$_H$-cover number.
\end{prop}

\begin{proof}[Sketch of the proof]
	The proof  is a straightforward generalisation of the LP-duality between the fractional matching and the fractional vertex-cover. Assume~$G$ is an $n$-vertex graph and assume the number of homomorphisms $h\in G^H$ of~$H$ in~$G$ is~$m$ (we have $m\le n^{|V(H)|}$). To attain the  fractional hom$_H$-tiling number is equivalent to the following instance of linear programming:
	\begin{align*}
	\mathrm{maximise}\quad&\quad\mathbf{v}^T\mathbf{z}\\
	\mathrm{subject~to}\quad&\quad A\mathbf{z}\le \mathbf{u}\\
	\mathrm{and}\quad&\quad \mathbf{z}\ge \mathbf{0}\;,
	\end{align*}
	where $\mathbf{v}$ and $\mathbf{u}$ are all-one vectors of length~$m$ and~$n$, respectively, $\mathbf{z}$ is the vector (of length~$m$) of variables (to be determined) corresponding to a fractional hom$_H$-tiling, and~$A$ is an $n\times m$ matrix. Each entry of the vector~$\mathbf z$ corresponds to the weight given to the  associated homomorphism.  Each column of~$A$ corresponds to a homomorphism $h\in G^H$ with image~$H'=h(H)$, each row corresponds to a vertex~$v$ of $G$, and $A_{H',v}=|h^{-1}(v)|$. Since this optimisation problem clearly has an optimal solution, the Strong LP-duality theorem says that it is equivalent to the following dual problem:
	\begin{align*}
	\mathrm{minimise}\quad&\quad\mathbf{u}^T\mathbf{y}\\
	\mathrm{subject~to}\quad&\quad A^T\mathbf{y}\ge \mathbf{v}\\
	\mathrm{and}\quad&\quad \mathbf{y}\ge \mathbf{0}\;,
	\end{align*}
	where  $\mathbf{y}$ is the vector (of length~$n$) of variables (to be determined) corresponding to a fractional hom$_H$-cover.
	\end{proof}

The core of the proof of Theorem~\ref{thm:result} is the next proposition, which  shows that a graph satisfying conditions similar to those in Theorem~\ref{thm:result} has large fractional cover. 

\begin{prop}\label{prop:r-colorable-cover}
Let $H$ be an $r$-colourable graph with colour class sizes \begin{equation}\label{eq:sizeOFell}\ell_1\ge \ldots\ge \ell_r>0\;\end{equation} and let $x\in (0,\frac{1}{|V(H)|})$. Set 
\begin{equation}\label{eq:def-delta} \delta:= (r-2+x\ell_r)\frac{n}{r-1}\;.\end{equation}

Then any $n$-vertex graph~$G$ with at least $(r-2+x|V(H)|)\frac{n}{r-1}$ vertices of degree at least~$\delta$ has fractional hom$_H$-cover number least~$xn$.
\end{prop}

\section{Proof of Theorem~\ref{thm:result}}\label{sec:proof_result}
Before jumping into the proof of Theorem~\ref{thm:result}, we explain first its main ideas.
By applying Szemer\'edi's regularity lemma to the input graph~$G$, we obtain a corresponding cluster graph~$\mathbf G$, which satisfies degree conditions similar to the ones of~$G$. Proposition~\ref{prop:r-colorable-cover} implies that the cluster graph~$\mathbf G$ has large  fractional hom$_H$-cover number. By the LP-duality for hom$_H$, the cluster graph~$\mathbf G$ has large fractional hom$_H$-tiling number. The embedding lemma yields a large $H$-tiling in the original input graph~$G$.

Let $H$ be a fixed $r$-colourable graph with colour classes $\ell_1\ge \cdots\ge \ell_r$ and let $x\in (0,\frac{1}{|V(H)|})$ be fixed. Assume we are given $\eta >0$. Notice that we can assume that $\eta <1$. Set $d:= \frac{1}{8}\left(  \frac{\eta x\ell_r}{r-1}\right)^2$. Since $\eta <1$ and $x<\frac{1}{|V(H)|}<\frac{1}{\ell_r}$, we have that $d<1$. Lemma~\ref{lem:emb} with input $d_{L\ref{lem:emb}}:= d/2$ and $\Delta_{L\ref{lem:emb}}:= |V(H)|$ outputs an $\varepsilon_0>0$. We define

\begin{equation}\label{eq:choice-eps}
\varepsilon_R:= \min\left\{\varepsilon_0, \frac{d}{6},\left(\frac{\eta}{6}\right)^2,\frac{1}{2}\left( \frac{1}{3|V(H)|x}-\frac{1}{3}\right)^2\right\}\;.
\end{equation}

Notice that $\left( \frac{1}{3|V(H)|x}-\frac{1}{3}\right)>0$, which implies that $\sqrt{\varepsilon_R}< \frac{1}{3|V(H)|x}-\frac{1}{3}$. 

The regularity lemma (Lemma~\ref{lem:RL}) with input $\varepsilon_R$
and $\ell:=\frac{1}{\varepsilon_R}$ outputs  $n_R,M\in \mathbb N$. Set

\begin{equation}\label{eq:choice-n}
 n_0:=\max\left\{n_R, \frac{2^{\Delta_{L\ref{lem:emb}}+1}\cdot|V(H)|\cdot M}{d^{\Delta_{L\ref{lem:emb}}}\cdot\sqrt{\varepsilon_R}\cdot(1-\varepsilon_R)}\right\}\;.
\end{equation}

Let $n\ge n_0$ and let~$G$ be a $n$-vertex graph with at least $(1+\eta)(r-2+x|V(H)|)\frac{n}{r-1}$ vertices of degree at least $(1+\eta)\delta$. We apply Lemma~\ref{lem:RL} on~$G$ (with parameters $\varepsilon_R$ and $\ell$) and obtain an equitable $\varepsilon_R$-regular partition $U_0, U_1, \ldots, U_p$ with $\ell<p<M$. 
We erase all edges within clusters~$U_i$, in irregular pairs, in regular pairs of density smaller than~$d$, and edges incident to the cluster~$U_0$.
Slightly abusing notation, we still call this subgraph~$G$. Let~$\mathbf G$ be the corresponding cluster graph with parameters $\varepsilon_R$, $m:=|U_i|$ ($i=1,\ldots,p$), and $d$.

\begin{cl}
After erasing edges as described in the previous paragraph, $G$ has at least~$(1+\eta/2)(r-2+x|V(H)|)\frac{n}{r-1}$ vertices of degree at least~$(1+\eta/2)\delta$.	
\end{cl}

\begin{proof}
By erasing all edges within clusters, we remove at most $\binom{m}{2}p$ edges. Notice that

\begin{equation*}
\binom{m}{2}p \leq \binom{n/p}{2}p \leq \frac{n^2}{2p}\leq \frac{\epsilon_Rn^2}{2}\,.
\end{equation*}

By erasing all edges between irregular pairs, we remove at most $\epsilon_Rp^2m^2\leq \epsilon_R n^2$ edges.
By erasing all edges between regular pairs of density smaller than $d$, we remove at most $\binom{p}{2}dm^2\leq \frac{d}{2}n^2$ edges.
Finally, by erasing all edges incident to cluster $U_0$, we remove at most $\epsilon_R n^2$ edges.
In total, the number of edges that have been removed from $G$ is at most 
\begin{equation*}
3\epsilon_R n^2+\frac{d}{2}n^2\leq dn^2= \frac{1}{8}\left(  \frac{\eta x\ell_r}{r-1}\right)^2n^2\,.
\end{equation*}

If the statement of the claim were not true, after the transformations, more than $\frac{\eta}{2}(r-2+x|V(H)|)\frac{n}{r-1}$ vertices of $G$ would have had their degree decreased by more than $\frac{\eta}{2}\delta$. This would imply that the number of edges removed from $G$ is at least

\begin{equation*}
\frac{1}{2}\frac{\eta}{2}(r-2+x|V(H)|)\frac{n}{r-1}\,\frac{\eta}{2}(r-2+x\ell_r)\frac{n}{r-1}> \frac{1}{8}\left(  \frac{\eta x\ell_r}{r-1}\right)^2n^2\,,
\end{equation*}
which yields a contradiction.
\end{proof}

\begin{cl} \label{cl:deg-cluster}
The cluster graph~$\mathbf G$ has at least~$(1+\eta/2)(r-2+x|V(H)|)\frac{p}{r-1}$ vertices of degree at least~$(1+\eta/2)(r-2+x\ell_r)\frac{p}{r-1}$.	
\end{cl}

\begin{proof}
For every vertex $v$ of $G$, we denote by $U_v$ the cluster of the regular partition containing~$v$. 

Let~$v$ a vertex of $G$ with degree at least~$(1+\eta/2)\delta$. Given a distinct vertex~$w$ of~$G$, we have that $vw\in E(G)$ implies that $U_vU_w\in E(\mathbf G)$. Since~$v$ has at most~$m$ neighbours in every cluster~$U_w$, we obtain that $U_v\in V(\mathbf G)$ has degree at least

\begin{equation*}
(1+\eta/2)\frac{\delta}{m}\geq (1+\eta/2)(r-2+x\ell_r)\frac{p}{r-1}\,.
\end{equation*}  

Similarly, let~$U_v\in V(\mathbf{G})$ have degree at least~$(1+\eta/2)(r-2+x\ell_r)\frac{p}{r-1}$. At most~$m$ vertices of~$G$ of degree at least~$(1+\eta/2)\delta$ belong to cluster~$U_v$. Therefore, the number of vertices of~$\mathbf G$ with degree at least~$(1+\eta/2)(r-2+x\ell_r)\frac{p}{r-1}$ is at least

\begin{equation*}
\frac{1}{m}(1+\eta/2)(r-2+x|V(H)|)\frac{n}{r-1}\geq (1+\eta/2)(r-2+x|V(H)|)\frac{p}{r-1} \,.
\end{equation*} 
\end{proof}

\begin{cl}\label{cl:tilingclustertograph}If the fractional hom$_H$-tiling number of~$\mathbf G$ is at least~$(1+3\sqrt{\varepsilon_R})xp$, then the $H$-tiling number of~$G$ is at least~$xn$. 
\end{cl}	

\begin{proof}
Let $f$ be a fractional hom$_H$-tiling of~$\mathbf G$ of maximum size. By hypothesis, \begin{equation*}
\sum f(h(H))\geq \left(1+3\sqrt{\varepsilon_R} \right)xp,
\end{equation*} 
where the sum runs through all homomorphisms $h\in \mathbf G^H$. 

Using the embedding lemma (Lemma~\ref{lem:emb}), we will sequentially find copies of~$H$ in~$G$, as long as $\sqrt{\varepsilon_R}m$ unused vertices are left in each cluster. We denote by~$W$ the vertices of~$G$ used by the copies of~$H$. At the beginning $W=\emptyset$. We define an auxiliary graph ${G}'\subseteq G$ as follows. Let $m'=\min\{|U_i\setminus W|\::\: i=1, \ldots, p\}$ be the smallest number of free vertices in a cluster. For each $i=1, \ldots, p$ let $U_i'$ be an arbitrary set of size $m'$ in $U_i\setminus W$. The graph~$G'$ is the subgraph of~$G$ induced by $\bigcup_{i=1}^pU_i'$. Observe that as long as $m'\ge \sqrt{\varepsilon_R}m$, we have by Lemma~\ref{lem:subdivide} that the pairs $(U_i',U_j')$ are $\varepsilon'$-regular with density at least $d'$, where $\varepsilon'=\sqrt{\varepsilon_R}$ and $d'=d/2$.

For every $h\in \mathbf G^H$, we shall construct $mf(h(H))(1-\sqrt{\varepsilon_R})$ vertex disjoint copies  of~$H$ in~$G$. Any such copy $g:H\rightarrow G$ will satisfy that, for every vertex~$v\in V(H)$,  its image $g(v)$ is contained in the cluster of~$G$ corresponding to~$h(v)$. So assume that we have a homomorphism $h\in \mathbf G^H$ and let $h(H)\subseteq \mathbf G$. We first construct a copy of~$H$ in~$\mathbf{G}_{|V(H)|}$ as follows: Let $U_i\in V(h(H))$; then each of the vertices in $h^{-1}(U_i)$ is embedded to a distinct vertex in the clone-set $U_i^{(s)}$ of $\mathbf{G}_{|V(H)|}$. 
Since $h$ is a homomorphism, no two vertices of this group are adjacent in~$H$. By repeating this procedure for all vertices in~$V(h(H))$, we obtain an embedding of the vertices of $H$ in $\mathbf{G}_{|V(H)|}$. By construction of $\mathbf{G}_{|V(H)|}$ and the fact that homomorphisms map edges of~$H$ to edges of~$h(H)$, we obtain that the edges of $H$ are also correctly embedded in $\mathbf{G}_{|V(H)|}$.
Due to the choices of $\varepsilon_R$ and $n_0$ given in equations~(\ref{eq:choice-eps}) and~(\ref{eq:choice-n}), we can use Lemma~\ref{lem:emb} with $d_{L\ref{lem:emb}}=d'$, $\Delta_{L\ref{lem:emb}}:=|V(H)|$, $m_{L\ref{lem:emb}}:= m'$, and $s_{L\ref{lem:emb}}:=|V(H)|$ to find an embedding of~$H$ in~$G$ such that for $U_i\in V(\mathbf G)$ all of the vertices in $h^{-1}(U_i)$ are embedded in the cluster $U_i$ of~$G$.
We put the vertices involved in this copy of $H$ in $G$ in the set~$W$. 
We repeat the process above $mf(h(H))(1-\sqrt{\varepsilon_R})$ times, and obtain $mf(h(H))(1-\sqrt{\varepsilon_R})$ vertex-disjoint copies of $H$ in $G$. This can be done using Lemma~\ref{lem:emb} (with the same parameters as above) as long as $m'\ge \sqrt{\varepsilon_R}m$. 

We next repeat this process for the remaining  homomorphisms  $h\in \mathbf G^H$. For this to be possible, we need to argue that $m'\ge \sqrt{\varepsilon_R}m$.
Let~$U_i$ be a vertex in $V(\mathbf G)$. The number of vertices in $U_i\setminus U'_i$ is at most the number of vertices used in copies of $H$ in $G$ so far. This is always at most $\sum mf(h(H))(1-\sqrt{\varepsilon_R})|(h)^{-1}(v)|$,
where the sum runs through all homomorphisms  $h\in \mathbf G^H$. In the next formula, the sums always run through all homomorphisms $h\in \mathbf G^H$:

\begin{align*} 
\sum mf(h(H))(1-\sqrt{\varepsilon_R})|(h)^{-1}(v)| ~& =~ m(1-\sqrt{\varepsilon_R})\left(\sum f(h(H))|(h)^{-1}(v)|\right) \\ 
& \leq ~  m(1-\sqrt{\varepsilon_R})\,. 
\end{align*}
Consequently, the number of vertices in $U'_i$ is at least
$
m-m(1-\sqrt{\varepsilon_R})= \sqrt{\varepsilon_R}m$.

It only remains to bound the size of the $H$-tiling of $G$ obtained by the above procedure. As we have seen, each homomorphism $h\in \mathbf G^H$ yields  $mf(h(H))(1-\sqrt{\varepsilon_R})$ vertex-disjoint copies of~$H$ in~$G$. Hence, the total size of the tiling is

\begin{equation*}
\sum mf(h(H))(1-\sqrt{\varepsilon_R})\geq m\left(1+3\sqrt{\varepsilon_R} \right)xp\left(1-\sqrt{\varepsilon_R}\right)\geq xn\left(1-\varepsilon_{R}\right)\left(1+3\sqrt{\varepsilon_R} \right)\left(1-\sqrt{\varepsilon_R}\right)\geq xn,
\end{equation*} 
where the sum runs through all homomorphisms $h\in \mathbf G^H$, and where we have used that $\varepsilon_R\leq 1/5$. Notice that this holds because $\varepsilon_R\leq d/6 <1/6$.
\end{proof}

From Proposition~\ref{prop:hom-duality}, we infer that in order to prove Theorem~\ref{thm:result} it is enough to show that the fractional hom$_H$-cover number of~$\mathbf G$ is at least~$(1+3\sqrt{\varepsilon_R})xp$. 
We apply Proposition~\ref{prop:r-colorable-cover} to the cluster graph~$\mathbf G$, where the value $x_{P\ref{prop:r-colorable-cover}}$ used in the proposition is $x_{P\ref{prop:r-colorable-cover}}:=(1+3\sqrt{\varepsilon_R})x$. For this to be valid, we need to argue that $(1+3\sqrt{\varepsilon_R})x\in (0,\frac{1}{|V(H)|})$, and that~$\mathbf G$ has the appropriate number of vertices of the appropriate degree.
Since $\varepsilon_R<\left( \frac{1}{3|V(H)|x}-\frac{1}{3}\right)^2$ (see~(\ref{eq:choice-eps}) and the paragraph afterwards), we have that $(1+3\sqrt{\varepsilon_R})x<\frac{1}{|V(H)|}$. Regarding the degree condition, by Claim~\ref{cl:deg-cluster}, $\mathbf G$ has at least $$(1+\eta/2)(r-2+x|V(H)|)\frac{p}{r-1}\overset{\eqref{eq:choice-eps}}{\ge} (1+3\sqrt{\varepsilon_R})(r-2+x|V(H)|)\frac{p}{r-1}$$ vertices of degree at least $$(1+\eta/2)(r-2+x\ell_r)\frac{p}{r-1}\overset{\eqref{eq:choice-eps}}{\ge} (1+3\sqrt{\varepsilon_R})(r-2+x\ell_r)\frac{p}{r-1}\;.$$ 

Thus, we can use Proposition~\ref{prop:r-colorable-cover} with the desired parameters, which ensures a fractional  hom$_H$-cover number of at least $(1+3\sqrt{\varepsilon_R})n$ in~$\mathbf G$. Proposition~\ref{prop:hom-duality} (LP-duality for hom$_H$) then implies that there is a fractional hom$_H$-tiling of size at least $(1+3\sqrt{\varepsilon_R})n$ in~$\mathbf G$. Finally, Claim~\ref{cl:tilingclustertograph} ensures an $H$-tiling in~$G$ of size at least~$xn$.

\section{Proof of Proposition~\ref{prop:r-colorable-cover}}\label{sec:proof_cover}

Let $c$ be any fractional hom$_H$-cover of $G$, and let $\left\| c \right\|$ denote its size. In order to prove the proposition, we want to show that $||c||\ge xn$.
 In the proof, we show that there is a clique of size~$r$ in~$G$, which is a homomorphic image of~$H$. This implies a lower bound on the total cover of the vertices of the clique. Among all possible such cliques, we consider one with minimum possible cover on its vertices with respect to the lexicographical order. From this we derive a lower bound for the cover of particular subsets of~$V(G)$, allowing us to compute a lower-bound for~$\|c\|$.

Let $L$ denote the set of vertices of $G$ having degree at least~$\delta$, and let $S$ denote the set containing the remaining vertices of $G$. By hypothesis,

\begin{equation}\label{eq:size-L}
\left|L \right|\geq \left(r-2+x|V(H)|\right)\frac{n}{r-1} \, .
\end{equation}

For any vertex $v$ of $G$, we denote the set of neighbours of $v$ in $L$ and $S$ by $N_L(v)$ and $N_S(v)$, respectively.

\begin{cl}\label{cl:existence_clique_in_L}
There is a clique~$K_{r-1}$ of size $r-1$ in $L$ and the vertices of any clique of size~$r-1$ in~$L$ have a non-empty common neighbourhood in $V(G)$. 
\end{cl}

\begin{proof}
We  find vertices $v_1,v_2,\ldots,v_{r-1}\in L$  of the clique $K_{r-1}$ in the following way. As~$L$ is non-empty, we pick~$v_1$ arbitrarily in $L$. If $r=2$, $v_1=v_{r-1}$ and we are done. If $r>2$,  suppose that we have already selected $i<r-1$ vertices $v_1,\ldots,v_i$ in $L$.
We define:
\begin{align*}
   N_L(i) ~& =~ \bigcap_{j=1}^{i} N_L(v_j),\\
	 N_S(i) ~& =~ \bigcap_{j=1}^{i} N_S(v_j)\, .
\end{align*}
We show that $|N_L(i)|>0$, as long as $i<r-1$. It is well-known that, given a set~$U$ of size~$m$ and~$i$ subsets~$A_j$ of~$U$,
\begin{align*}
\left|\bigcap_{j=1}^{i} A_j \right| ~ \geq ~ \sum_{j=1}^{i}\left|A_j\right|-(i-1)\cdot m \, .
\end{align*}
Consequently,
\begin{align} 
\left|N_L(i)\right|+\left|N_S(i)\right| ~& \geq~ i \cdot \delta -(i-1)\cdot n=i \cdot \left( (r-2+x\ell_r)\frac{n}{r-1}\right) -(i-1)\cdot n \nonumber\\ \nonumber
& = ~ \left(\frac{r-1-i}{r-1}\right)\cdot n +\frac{i x \ell_r n}{r-1}\\ 
& = ~ \left(r-1-i+ i x \ell_r\right)\frac{n}{r-1}  \,. \label{eq:inc-exc}
\end{align}
Note that Inequality~(\ref{eq:inc-exc}) holds also in the case when $i=r-1$.

Since
\begin{align*}
\left|N_S(i)\right| ~ \leq ~ \left|S\right|=n- \left|L\right|\overset{\eqref{eq:size-L}}{\leq} (1-x|V(H)|)\frac{n}{r-1}
\end{align*}
and $i<r-1$, we have that $\left|N_L(i)\right|>0$ and we can pick $v_{i+1}\in N_L(i)$. Continuing in this way, we obtain $K_{r-1}$ in~$L$.

Observe that if $i=r-1$, we obtain that $v_1, \ldots, v_{r-1}$ have common neighbourhood of size at least $x\ell_r n>0$  (and hence there is a clique of size~$r$ in~$G$). However, this common neighbourhood may lie completely in~$S$.
\end{proof}

Among all possible cliques of size $r-1$ in $L$, pick one with the smallest cover with respect to the lexicographical order and let $v_1, \ldots, v_{r-1}$ be its vertices (in lexicographical order w.r.t.\ the cover). We set $u\in  N_L(r-1)$ to be a vertex with the smallest $c(u)$, and $w\in N_S(r-1)$
 one with the smallest~$c(w)$. 
   It may happen that~$u$ or~$w$ does not exist, but at least one of them does. For $i=1,\ldots,r-1$, we define $\alpha_i=c(v_{i})$. Let $\alpha_L=c(u)$, if $u$ exists, and set $\alpha_L=1$, otherwise. Similarly, let $\alpha_S=c(w)$, if $w$ exists, and set $\alpha_S=1$, otherwise. Let $\alpha_r= \min\{\alpha_L,\alpha_S\}$. 
   We then have:

\begin{equation}\label{eq:alphas}
\alpha_1\leq \alpha_2\leq\cdots\leq \alpha_{r-1}\leq \alpha_L\,.
\end{equation}

As~$H$ is $r$-colourable, there exists a homomorphism $h:H\rightarrow G$ that sends all vertices in the $\ell_i$-sized color class of $H$ to $v_i$ (for $i=1,\ldots,r-1$), and all vertices in the $\ell_r$-sized color class of $H$ to $u$ or $w$. Since $c$ is a fractional hom$_H$-cover of $G$, we obtain the following relation:
\begin{equation}\label{eq:cover-r}
\ell_1\alpha_1+\ell_2\alpha_2+\cdots+\ell_r\alpha_r\ge 1\,.
\end{equation}

As $x<\frac{1}{|V(H)|}$, we have 
\begin{equation}\label{eq:size-n}
n> x\left|V(H)\right|n\,.
\end{equation}

To obtain a lower bound for $\| c \|$, we will use the following auxiliary inequality:

\begin{cl}\label{cl:caseAnalisis}
	$\left|   N_L(r-1)   \right| (\alpha_L-\alpha_{r-1}) + 
	\left| 	N_S(r-1)	\right| \cdot \alpha_r \ge \alpha_r\cdot x\cdot \ell_r \cdot n -\alpha_{r-1}\left[	(r-1)\delta -(r-2)n	\right]
	$.
\end{cl}

\begin{proof}
	In order to prove the claim, we need to consider the two following cases:
	\begin{itemize}
		\item[CASE 1:] $\alpha_L\ge  \alpha_r+\alpha_{r-1}$
		\item[CASE 2:] $\alpha_L< \alpha_r+\alpha_{r-1}$
	\end{itemize}
	In CASE 1, we have 
	\begin{align}
		\left|   N_L(r-1)   \right| (\alpha_L-\alpha_{r-1}) + 
		\left| 	N_S(r-1)	\right| \cdot \alpha_r  & \ge  \left| N_L(r-1) \right|\cdot\alpha_r + \left|  N_S(r-1) \right| \cdot  \alpha_r \nonumber
	\\	& \hspace*{-0.37cc} \overset{\eqref{eq:inc-exc}}{\ge} \alpha_r\cdot x \cdot \ell_r\cdot n\nonumber\\
	& \ge \alpha_r\cdot x \cdot \ell_r\cdot n +\alpha_{r-1}\left( (r-2)n-(r-1)\delta     \right)\,. \nonumber
	\end{align}
	
	In CASE 2, we get 
	\begin{align}
	\MoveEqLeft[4] \left|   N_L(r-1)   \right| (\alpha_L-\alpha_{r-1}) + 
		\left| 	N_S(r-1)	\right| \cdot \alpha_r\nonumber\\
		& \ge   \left|   N_L(r-1)   \right| (\alpha_L-\alpha_{r-1}) + 
		\left| 	N_S(r-1)	\right| (\alpha_L-\alpha_{r-1})\nonumber\\
	& \hspace*{-0.37cc} \overset{\eqref{eq:inc-exc}}{\ge}  
\left((r-1)\delta-(r-2)n\right)(\alpha_L-\alpha_{r-1})\nonumber	\\
&=\alpha_{r-1}\left(  (r-2)n-(r-1)\delta  \right)+ \alpha_L\cdot x \cdot \ell_r\cdot n\nonumber\\
&	\ge \alpha_{r-1}\left(  (r-2)n-(r-1)\delta  \right)+
	\alpha_r\cdot x \cdot \ell_r\cdot n\,. \nonumber
	\end{align}

\end{proof}

We now have the necessary ingredients to prove that $||c||\ge xn$. We lower-bound $\|c\|$ by computing the cover of particular neighbourhoods of $v_1,\ldots,v_{r-1}$.

\begin{align}
 \left\| c \right\| = & \sum_{v\in V(G)}c(v)\nonumber\\
\geq&~\left| L \right| \cdot \alpha_1 + \sum_{i=1}^{r-2} \left| N_L(i)\right| \cdot (\alpha_{i+1}- \alpha_i)+       \left| N_L(r-1)\right| \cdot (\alpha_{L}- \alpha_{r-1})+
       \left|  N_S(r-1) \right| \cdot \alpha_S\,.
\label{eq:constantCover}
\end{align}

By~\eqref{eq:alphas}, we have that $(\alpha_{i+1}-\alpha_i)\ge 0$, for $i=1,\ldots,r-2$ and that $\alpha_L-\alpha_{r-1}\ge 0$. Therefore,
\begin{align}
\eqref{eq:constantCover}\overset{\eqref{eq:inc-exc}}\ge &\left| L \right| \cdot \alpha_1 + \sum_{i=1}^{r-2} \left[ i \cdot \delta -(i-1)\cdot n - \left| N_S(i) \right| \right] \cdot (\alpha_{i+1}- \alpha_i)\nonumber \\
& + \left| N_L(r-1)\right| \cdot (\alpha_{L}- \alpha_{r-1})+
\left|  N_S(r-1) \right| \cdot \alpha_S\nonumber \\
\geq & ~ \left| L \right| \cdot \alpha_1 + \sum_{i=1}^{r-2} \left[ i \cdot \delta -(i-1)\cdot n - \left| S \right| \right] \cdot (\alpha_{i+1}- \alpha_i)\nonumber \\
&   +   \left| N_L(r-1)\right| \cdot (\alpha_{L}- \alpha_{r-1})+
\left|  N_S(r-1) \right| \cdot  \alpha_r\nonumber\\
=&~   \alpha_1(|L|+|S|-\delta)+\sum_{i=2}^{r-2} \left(  \alpha_i\cdot \left(n-\delta \right)   \right)+\alpha_{r-1}\left(n-(r-2)(n-\delta)-|S|   \right)\nonumber \\
&~+  \left| N_L(r-1)\right| \cdot (\alpha_{L}- \alpha_{r-1})+
\left|  N_S(r-1) \right| \cdot  \alpha_r\,. \label{eq:replacebyS}
\end{align}

Using the auxiliary inequality from Claim~\ref{cl:caseAnalisis}, we get
\begin{align}
\eqref{eq:replacebyS}\ge ~& ~\alpha_1 \cdot \left(\left| L \right| + \left|S\right| - \delta \right)+\sum_{i=2}^{r-2} \left(  \alpha_i\cdot \left(n-\delta \right)   \right)+\alpha_{r-1}\cdot \left(n-\delta - \left|S\right|\right) + \alpha_r\cdot x\cdot \ell_r \cdot n	\nonumber	\\
					\overset{\eqref{eq:cover-r}}{\geq}& ~\alpha_1 \cdot \left(n - \delta \right)+\sum_{i=2}^{r-2} \left(  \alpha_i\cdot \left(n-\delta \right)   \right)+\alpha_{r-1}\cdot \left(\left|L\right|-\delta \right) + \frac{1}{\ell_r}\cdot \left(1- \sum_{i=1}^{r-1} \ell_i \cdot \alpha_i\right) \cdot x\cdot \ell_r \cdot n\nonumber \\
					 \overset{\eqref{eq:def-delta},\eqref{eq:size-L}}\ge  &~ \sum_{i=1}^{r-2} \left(  \alpha_i\cdot \left(n-\delta - x \cdot\ell_i\cdot  n   \right)   \right)+\alpha_{r-1}\cdot \left(x \left(\left|V(H)\right| - \ell_r\right)\frac{n}{r-1}  -x\cdot \ell_{r-1}\cdot n\right)+x \cdot n \nonumber \\
  \overset{\eqref{eq:def-delta}}{=} & \sum_{i=1}^{r-2} \left(  \alpha_i\cdot \left( (1-x\cdot \ell_r)\frac{n}{r-1} - x \cdot\ell_i\cdot  n   \right)   \right)+\alpha_{r-1}\cdot \left(x \left(\left|V(H)\right| - \ell_r\right)\frac{n}{r-1}  -x\cdot \ell_{r-1}\cdot n\right)+x \cdot n\nonumber \\
 \overset{\eqref{eq:size-n}}{>} &~ x \cdot n+ \sum_{i=1}^{r-1}\left( \frac{\alpha_i}{r-1} \left(\sum_{j=1}^{r-1}\left(x\cdot \ell_j \cdot n\right)-(r-1)\cdot x\cdot \ell_i\cdot n  \right) \right)\nonumber\\
 = ~& ~ x \cdot n+ \frac{x \cdot n}{r-1}\cdot \sum_{i=1}^{r-1}\left( \alpha_i  \cdot \left(\sum_{j=1}^{r-1}\left(\ell_j - \ell_i  \right) \right)\right)\nonumber\\	 = ~& ~ x \cdot n+ \frac{x \cdot n}{r-1}\cdot \sum_{i=1}^{r-1} \sum_{j=1}^{i-1}\left(\left(\ell_j - \ell_i  \right) \cdot \left(\alpha_i - \alpha_j  \right)\right)\nonumber \\
 \overset{\eqref{eq:alphas}\, \eqref{eq:sizeOFell}}{\geq} & ~ x \cdot n \label{eq:atleastk}\, .\nonumber
    \end{align}

Hence the fractional hom$_H$-cover number of $G$ is at least $xn$.

\section{Final remarks}

The approximation in Theorem~\ref{thm:result} cannot be totally removed, as witnessed by the following graphs. Assume that the graph $H=K_{3,3,3}$ we want to tile is the balanced complete  tripartite graph on nine vertices.
The vertex set of the host graph~$G$ is partitioned into four sets $V_1, V_2,V_3,V_4$, with $|V_1|=3xn-2$, $|V_2|=3xn+2$, $|V_3|=(1-3x)n/2$ and $|V_4|=(1-9x)n/2$, where the sets $V_1,V_3$, and~$V_4$ are independent, the graph induced by~$V_2$ is a (spanning) cycle, and all the edges between~$V_3$ and $V_1\cup V_2\cup V_4$, and between~$V_1$ and~$V_2$ are present. The graph~$G$
has $(r-2+x|V(H)|)\frac{n}{r-1}$ vertices of degree at least $\delta:=(r-2+x\ell_r)\frac{n}{r-1}$ but does not contain~$xn$ vertex disjoint copies of~$K_{3,3,3}$, for $x>1/n$. 

It would be interesting however to know what approximation is really needed, i.e., given a graph~$H$ and~$x\in (0,\frac{1}{|V(H)|})$, what functions~$f_1$ and~$f_2$ guarantee that any $n$-vertex graph with at least $\left(r-2+x\frac{1}{|V(H)|}\right)\frac{n}{r-1}+f_1(n)$ vertices of degree at least $\delta+f_2(n)$ contains an $H$-tiling of size at least~$xn$? In particular, do $f_i(n)$, $i=1,2$ depend on~$n$?

\section{Acknowlegement}

We thank an anonymous referee for his advice on how to improve the presentation of the paper. We also thank J.~Hladk\'y for useful discussions.

\end{document}